\documentclass[12pt]{amsart}
\usepackage{amsfonts}

\usepackage[letterpaper]{geometry}

\usepackage{amssymb}
\def\N{\mathbb{N}}

\newtheorem{thm}{\bf Theorem}[section]
\newtheorem{lemma}[thm]{\bf Lemma}

\newtheorem{cor}[thm]{\bf Corollary}
\newtheorem{st}[thm]{\bf Statement}

\theoremstyle{definition}

\newtheorem{que}{Question}

\begin{document}

\title{An elementary question of Erd\H os and Graham}

 \author[Norbert Hegyv\'ari]{Norbert Hegyv\'ari}
 \address{Norbert Hegyv\'{a}ri, ELTE TTK,
E\"otv\"os University, Institute of Mathematics, H-1117
P\'{a}zm\'{a}ny st. 1/c, Budapest, Hungary and associated member of Alfr\'ed R\'enyi Institute of Mathematics, Hungarian Academy of Science, H-1364 Budapest, P.O.Box 127.}
 \email{hegyvari@renyi.hu}

\begin{abstract} Let $A_k=\{r(k-r): 1\leq r \leq k-1\}$. Erd\H os and Graham asked about the cardinality of the set of common elements. We answer this elementary question and apply our result to a sum-product type result.

MSC 2020 Primary   Secondary 

Keywords: Erd\H os question, sum-product type result
\end{abstract}

 \maketitle
\section{Introduction}

For an integer $k\geq 2$, let $A_k=\{r(k-r): 1\leq r \leq k-1\}$. 
In 1980 Erd\H os and Graham posed the following question ([1] and [2]:
\begin{que}
"Let $n,m$ be positive integers and consider the two sets $A_n,A_m$. Can one estimate the number of integers common to both? Is this number unbounded? It should certainly be less than $(mn)^\varepsilon$ for every $\varepsilon>0$ if $mn$ is sufficiently large."
\end{que}
In this short note we answer to this question in section 3. In the last section we apply the positive answer to a sum-product result (although this application gives a conditional result). The reason to show this result is that our method uses arithmetical argument and is close to the state-of-art result. The proof of  best results are proved using geometrical, and so called 'higher energy' method, and states for sets of real numbers. Our approach is elementary (and maybe the original question of Erd\H os and Graham was also related to the sum-product problem).

I think that, because of the 'rigidity' of the integers, the best lower bound on the range of integers will be orders of magnitude larger than the best lower bound in the real case.

\medskip

Both of the above questions are answered in the affirmative, in a sharper form. 

Let $T_{m,n}=\{(M,N): NM=m^2-n^2; \ M+N<2m; \ 0<M\leq N<M+2n\}$, and write $\tau_{m,n}:=|T_{m,n}|$.

\begin{thm}
Let $m>n$. Then $|A_n\cap A_m|\leq \tau_{m,n}$. Furthermore if $m$ is even and $n$ is odd then equality holds.
\end{thm}

\begin{cor}
For every $\varepsilon>0$ there exists a threshold $n_0$ such that $m>n>n_0$ we have 
$$
|A_n\cap A_m|<(m)^{\frac{(2+\varepsilon)\log 2}{\log\log m}}
$$
\end{cor}
\begin{thm}
For every $s\in \N$, there exists an infinite sequence $\{n_1<m_1<n_1<n_2<m_2<\dots\}$ such that 
$$
|A_{n_k}\cap A_{m_k}|=s; \quad  k=1,2,\dots
$$
\end{thm}

\section{Notation}
The cardinality of a finite set $X$ is denoted by $|X|$. $\N$ denotes the set of non-negative integers; $\N=\{0,1,2,\dots\}$. Denote by $d(x)$ the number of positive divisors of $x$; $d(x)=\sum_{d|x; x>0}1$.

\section{Proofs}

\begin{proof}[Proof of Theorem 1.1]
Let $m>n$. We only have to count the number of solutions of $k(n-k)=r(m-r)$ in the cases $1\leq k\leq n/2$ and $1\leq r\leq m/2$. This equation is equivalent to
$(m-2r)^2-(n-2k)^2=m^2-n^2$ or 
$$
(m-2r+n-2k)(m-2r-n+2k)=m^2-n^2.
$$
Let $M:=m-2r+n-2k$ and $N:=m-2r-n+2k$. It is easy to check that $N\geq M>0$. So, for a given factorization $MN=m^2-n^2$, we get
\begin{equation}\label{1}
r=\frac{2m-N-M}{4}; \quad k=\frac{2n+M-N}{4}.
\end{equation}
Furthermore we have $N+M=2m-2r<2m$ and $N<2n+M$. Hence we conclude if the pair $(k,r)$ is a solution of $k(n-k)=r(m-r)$ then for some $(M,N)\in T_{m,n}$, (\ref{1}) holds. Furthermore if $(k',r')$ is also a solution then the corresponding pair $(M',N')$ differs from $(M,N)$. 

This yields that there is an injection from the set $A_n\cap A_m$ to $T_{m,n}$, hence $|A_n\cap A_m|\leq |T_{m,n}|=\tau_{m,n}$.

To prove the second part of the theorem let $m=2s; \ n=2p+1$. We have $m^2-n^2=4(s^2-p^2-p)-1=4L-1$. Now let us take the pair $(M,N)$, for which $M|4L-1; \ N(4L-1)/M$, and fulfills the condition $M+N<2m; \ 0<M\leq N<M+2n$ (provided such a pair exists).

Since $MN\equiv-1 \pmod 4$ thus $M+N\equiv 0\pmod 4$, hence the numbers $k$ and $r$ in (\ref{1}) are integers.
    
\end{proof}

\begin{proof}[Proof of Corollary 1.2]
    It is well-known that for $x>x(\varepsilon)$, $d(x)<x^{\frac{(1+\varepsilon)\log 2}{\log\log x}}$ (see e.g. in [3] p.262, Thm 317). Since $\tau_{m,n}\leq d(m^2-n^2)$ the corollary follows.
\end{proof}
\begin{proof}[Proof of Theorem 1.3]
Let $s\in \N$ be a fixed integer.
Write $m=m_1$, $n=n_1$ and $p=p_1\geq 3$. Pick an $\alpha\in \N$, for which 
\begin{equation}\label{2}
    2^\alpha>p^s
\end{equation}
and let $m=2^\alpha p^s+1; \ n=m-2.$ Then $m^2-n^2=2^{\alpha+2}p^s$. Since $2m\equiv2n\equiv2\pmod 4$ by (\ref{1}) we infer that
\begin{equation}\label{3}
   N_i= 2^{\alpha+1}p^i; \quad M_i=2p^{s-i} \quad 0\leq i\leq s
\end{equation}
are the possible pairs.

We are going to show if $i\neq s$ then every pair of $(M_i,N_i)$ satisfies (\ref{3}) there is an unique pair of $(k,r)$ for which $k(n-k)=r(m-r)$. 
Taking into account (\ref{1}) we have to check that $M_i+N_i<2m$, which in our case is $2^{\alpha+1}p^i+2p^{s-i}<2^{\alpha+1}p^s+2$. This can also be written as $0<(2^{\alpha}p^i-1)(2p^{s-i}-1)$, which is clearly true.

Furthermore we check the inequalities $M_i\leq N_i<M_i+2n$. $M_i\leq N_i$ follows from ($\ref{2}$). 

For the second inequality we get $N_i-M_i<2^{\alpha+1}p^s-2=2n$. The equation $m^2-n^2=M_iN_i$ has $s$ many solution as we want.

Finally, for every pair $n_i<m_i; \ i=2,3\dots $ we choose a $p_i>p_{i-1}$, for which the previous process works and which guarantees that $A_{n_i}\neq A_{n_j}$ (and $A_{m_i}\neq A_{m_j}$) for $(j<i)$

\end{proof}

\section{The Erd\H os question and the sum-product problem}

One suspects that Erd\H os asked this question in order to get an estimate to the sum-product problem. In 1983 Erd\H os and Szemer\'edi gave a non-trivial estimation $\max\{|A+A|,|AA|\}>|A|^{1+c}; \ c>0$ where $A\subseteq \N$ is a finite set. 

Surprisingly, Elekes has taken a major step forward with a geometric tool using an incidence theorem. A further ingenious improvement was made by Solymosi (also by a "geometric aspect"), and finally Rudnev and Stevens have the 'state-of-art' result where the exponent is $4/3+0,0017\dots-o(1)$ (see [5],[6]). 

\subsection{Conditional sum-product problem}

It is expected that the above theorems are not strong enough to give a sum-product result. Nevertheless, for the sake of completeness, we show a conditional result that can be conclude from the result above.
\begin{st}
Let $0<c<1$ and $n>n(c)$ be an integer large enough. Let  $A$ be a set of integers with $n$ many elements and assume that $A\subset [1,n^{(\log\log n)^c}]$. Then 
$$
\max\{|A+A|,|AA|\}\gg n^{4/3-\frac{3}{(\log\log n)^{1-c}}}.
$$
\end{st}
\begin{proof}
Let $\sqcup_{i\geq 1}B_i$ be a partition of $A\times A$, where $B_k:=\{(a, a')\in A\times A: a+a'=k\}$. 

If $|A+A|\geq n^{4/3-\frac{3}{(\log\log n)^{1-c}}}$ then we are done. Assume to the contrary. Write $T=|A+A|$, so $T\leq n^{4/3-\frac{3}{(\log\log n)^{1-c}}}$. For $1\leq k\leq T$, let $U_k=\{B_k: |B_k|>0.5n^{2/3-\frac{2}{\log\log n}}\}$. Now 
$$
|(A\times A)\setminus \cup_{B_k\in U_k}B_k|<0.5n^{2/3-\frac{2}{\log\log n}}n^{4/3-\frac{2}{\log\log n}}<0.5n^2 
$$
so $|\cup_{B_k\in U_k}B_k|> 05n^2$. 

Thus we conclude that $m:=|U_k|\geq 0.5n$, otherwise 
$$
|\cup_{B_k\in U_k}B_k|\leq \sum_{B_k\in U_k}|B_k|\leq \sum_{B_k\in U_k}n<0.5n^2.
$$
\begin{lemma}
Let $\mathcal{H}=(V,E)$ be a hypergraph with $m$ many edges and assume that the cardinality of each edge is at least $R$. Assume that for every distinct edges $e,e'\in E$ $|e\cap e'|\leq k$. Then,  $|V|\geq \frac{mR^2}{R+(m-1)k}$.
\end{lemma}
It is an easy consequence of ([4, Chap. 13, Prob.13]).
\end{proof}
For every $B_k\in U_k$ let the corresponding product set be $C_k:=\{a\cdot a':(a,a')\in B_k\}$. Clearly if $a\cdot a'=a''\cdot a'''$ for some four-tuple from $C_k$ we have $\{a,a'\}=\{a'',a'''\}$.

Clearly $C_k\subseteq A_k$ and $C_k\cap C_j\subseteq A_k\cap A_j$. Since $A\subset [1,n^{(\log\log n)^c}]$ by Corollary 1.2 we get that
$$
|C_k\cap C_j|\ll (n)^{\frac{2}{(\log\log n^{c'}}} \quad k\neq j
$$
where c'=1-c and $C$ is an absolute constant.

Finally use Lemma 4.2 with $k=n^\frac{2}{(\log\log n)^{1-c}}$, $m:=|U_k|$, $R\geq 0.5n^{2/3-\frac{2}{\log\log n}}\}$. Observe, that $R\leq (m-1)k$, thus $V\geq \frac{mR^2}{2mk}=R^2/2k$, in which $|AA|\gg n^{4/3-\frac{3}{(\log\log n)^{1-c}}}$ if 
 $n$ is large enough, which is the desired lower bound.

\end{document}